\documentclass[article,10pt]{aiaa-pretty}    

\usepackage[utf8]{inputenc}  

\usepackage{amssymb}
\usepackage{amsmath}
\usepackage{booktabs}
 \usepackage{varioref}
 \usepackage{wrapfig}
 \usepackage{threeparttable}
 \usepackage{dcolumn}
  \newcolumntype{d}{D{.}{.}{-1}}
 \usepackage{nomencl}
  \makeglossary
 \usepackage{subfigure}
 \usepackage{graphicx}
 \usepackage{subfigmat}
 \usepackage{fancyvrb}
  \fvset{fontsize=\footnotesize,xleftmargin=2em}
 \usepackage{lettrine}
\usepackage{cite}

\newtheorem{theorem}{Theorem}
\newtheorem{proposition}{Proposition}
\newtheorem{remark}{Remark}

\newtheorem{definition}{Definition}

\newenvironment{proof}{{\it Proof. }}{\hfill $\Box$}

\newcommand{\Real}{\mathbb R}
\newcommand{\set}[1]{\left\{#1\right\}}
\newcommand{\real}[1]{{\mathbb R}^{#1}}
\newcommand{\bb}{{\boldsymbol b}}

\newcommand{\bC}{{\boldsymbol C}}

\newcommand{\be}{{\boldsymbol e}}
\newcommand{\bff}{{\boldsymbol f}}
\newcommand{\bg}{{\boldsymbol g}}
\newcommand{\bh}{{\boldsymbol h}}

\newcommand{\bp}{{\boldsymbol p}}

\newcommand{\bu}{{\boldsymbol u}}

\newcommand{\bv}{{\boldsymbol v}}
\newcommand{\bw}{{\boldsymbol w}}
\newcommand{\bx}{\boldsymbol x}

\newcommand{\bxf}{{\bx(\cdot)}}  

\newcommand{\buf}{{\bu(\cdot)}}  

\newcommand{\bA}{{\boldsymbol A}}
\newcommand{\bB}{{\boldsymbol B}}

\newcommand{\bP}{{\boldsymbol P}}

\newcommand{\U}{\mathbb{U}}

\newcommand{\bzero}{{\bf 0}}

\newcommand{\bX}{{\mbox{\boldmath $X$}}}
\newcommand{\bV}{{\mbox{\boldmath $V$}}}
\newcommand{\bU}{{\mbox{\boldmath $U$}}}

\newcommand{\bnu}{{\mbox{\boldmath $\nu$}}}
\newcommand{\bmu}{{\mbox{\boldmath $\mu$}}}

\newcommand{\blam}{{\mbox{\boldmath $\lambda$}}}

\author{ %
R. J. Proulx\thanks{Research Professor, Control and Optimization Laboratories, Space Systems Academic Group}
and
I. M. Ross\thanks{Distinguished Professor and Program Director, Control and Optimization, Department of Mechanical and Aerospace Engineering}
\\
\textit{Naval Postgraduate School, Monterey, CA 93943}
}
\title{Implementations of the Universal Birkhoff Theory for Fast Trajectory Optimization}

\abstract{This is part II of a two-part paper.  Part~I presented a universal Birkhoff theory for fast and accurate trajectory optimization. The theory rested on two main hypotheses.  In this paper, it is shown that if the computational grid is selected from any one of the Legendre and Chebyshev family of node points, be it Lobatto, Radau or Gauss, then, the resulting collection of trajectory optimization methods satisfy the hypotheses required for the universal Birkhoff theory to hold. All of these grid points can be generated at an $\mathcal{O}(1)$ computational speed.
Furthermore, all Birkhoff-generated solutions can be tested for optimality by a joint application of Pontryagin's- and Covector-Mapping Principles, where the latter was developed in Part~I.  More importantly, the optimality checks can be performed without resorting to an indirect method or even explicitly producing the full differential-algebraic boundary value problem that results from an application of Pontryagin's Principle. Numerical problems are solved to illustrate all these ideas.  The examples are chosen to particularly highlight three practically useful features of Birkhoff methods: (1) bang-bang optimal controls can be produced without suffering any Gibbs phenomenon, (2) discontinuous and even Dirac delta covector trajectories can be well approximated, and (3) extremal solutions over dense grids can be computed in a stable and efficient manner.
}

\begin{document}
\maketitle


\section{Introduction}\label{sec:Intro}
In \cite{newBirk-part-I}, a new theory was advanced for fast and accurate trajectory optimization.  This theory was based on the universal Birkhoff interpolants developed in \cite{newBirk-2023}. The universal Birkhoff theory of \cite{newBirk-part-I} rested on two hypotheses; namely, that there exists a grid such that:
\begin{enumerate}
\item the $a$- and $b$-versions of the Birkhoff interpolant are equivalent; and,
\item the Birkhoff quadrature formula is accurate to within any given small $\epsilon > 0$.
\end{enumerate}
In this paper, we first provide a proof of existence of such a grid; in fact, we show the existence of multiple grids that satisfy the two hypotheses. In addition, we provide implementation details for these particular choices of grid points and thereafter demonstrate the efficacy of the resulting Birkhoff methods by solving a set of illustrative problems.

As detailed in \cite{newBirk-part-I} and \cite{newBirk-2023}, the Birkhoff interpolants and the ensuing theory are different from the prior Birkhoff methods proposed in the literature\cite{wang,fastmesh,furtherResults}; in particular, it was shown in \cite{newBirk-part-I} that the universal interpolants designed in \cite{newBirk-2023} also satisfy the covector mapping principle (CMP)\cite{perspective,ross:roadmap-2005,ross-book,CMP:history,CMP-CDC-2006}, caveated by the two hypotheses noted earlier.  That is, if a trajectory optimization method based on the new Birkhoff interpolants proposed in \cite{newBirk-2023} is constructed, then, dualization and discretization can be commuted\cite{ross:roadmap-2005,ross-book,perspective}.  To explain, in simpler terms, why commutative operations in trajectory optimization are crucial\cite{ross:roadmap-2005} to generate the correct solution, consider the following simple univariate function,
\begin{equation}
g^N(x) := \frac{N x + x^2}{N^2}, \quad x \in \Real, \ N \in \mathbb{N}
\end{equation}
Suppose we minimize $g^N(x)$ with respect to $x$ (holding $N$ constant).  It is straightforward to show that this yields,
\begin{equation}\label{eq:mingN}
\min_x g^N(x) = -\frac{1}{4}
\end{equation}
Obviously, taking the limit of the left hand side of \eqref{eq:mingN} as $N \to \infty$ generates,
\begin{align}
 \lim_{N \to \infty}\left(\min_x g^N(x) \right) &= -\frac{1}{4}
\end{align}
%
That is, minimizing $g^N(x)$ before taking the limit $N \to \infty$ produces the number, $-1/4$; however, if we reverse the operations, we get a different answer:
\begin{subequations}
\begin{align}
 \lim_{N \to \infty} g^N(x) &= 0 \\
 \min_x \left( \lim_{N \to \infty} g^N(x) \right) &= 0
\end{align}
\end{subequations}
Obviously,
\begin{equation}
\lim_{N \to \infty}\min_x g^N(x) \ne \min_x \lim_{N \to \infty} g^N(x)
\end{equation}
That is, $\lim$ and $\min$ are not commutative operations.  Note that $\arg\min$ and $\lim$ are also not commutative:
\begin{equation}
\arg\min_x \left(\lim_{N \to \infty} g^N(x)\right) \ne \lim_{N \to \infty} \left(\arg\min_x g^N(x)\right)
\end{equation}
Trajectory optimization involves a number of commutative operations\cite{ross-book,ross:roadmap-2005,perspective}; hence, it is crucial that any proposed trajectory optimization method satisfy the CMP. Production of covector mapping theorems for the original Birkhoff method used in \cite{fastmesh,furtherResults} parallels that of previous pseudospectral methods\cite{LGL:jgcd-2001,acc:hybrid,advances,cheb-costate}; i.e.,  the development of intricate details that are specific to the different choices of grid points.  In a quest for a more universal theory akin to the unified Lagrange pseudospectral (PS) method\cite{arb-grid-aas,arb-grid,PSReview-ARC-2012}; i.e., methods based on differentiation matrices\cite{fornberg,boyd,trefethen-2000,hgg-2007},  it was shown in \cite{newBirk-part-I} that the new Birkhoff interpolants proposed in \cite{newBirk-2023} satisfy the CMP for any choice of grid points so long as two hypotheses are satisfied. With this perspective in mind, the main contributions of this paper are as follows:
\begin{enumerate}
\item A proof of existence of viable Birkhoff methods.  This involves proving that a family of Legendre and Chebyshev grids (Lobatto, Radau and Gauss) satisfy the two hypotheses proposed in \cite{newBirk-part-I}, and hence the CMP;
\item Implementation details of Birkhoff methods for the family of Legendre and Chebyshev grid points;
\item Results from select numerical problems to illustrate three features of the universal Birkhoff methods that are important in practical applications:
\begin{enumerate}
\item Production of bang-bang-type solutions without suffering any Gibbs phenomenon despite the computational technique being based on Legendre/Chebyshev polynomials;
\item Solution to state-constrained optimal control problems that require the existence and computation of discontinuous costates and Dirac-delta-type path covector trajectories; and
\item Generation of solutions over very dense grids 
without suffering round-off errors or extensive preconditioning.\label{point:N>1000}
\end{enumerate}
\end{enumerate}
The preceding point~\ref{point:N>1000} is particularly relevant when compared with PS methods based on Lagrange interpolants\cite{boyd,fornberg,trefethen-2000,hgg-2007}.  As noted in \cite{newBirk-part-I} and amplified in \cite{wang,newBirk-2023}, the main initial motivation for a Birkhoff theory was a need to produce a more perfect preconditioner\cite{hesthavan,elbarbary} to differentiation matrices.  The search for a perfect preconditioner led to an ``accidental'' new starting point for PS theory and trajectory optimization itself\cite{newBirk-part-I}.  Together with \cite{newBirk-part-I}, this paper provides a more complete picture of this emerging new theory.

\section{Legendre and Chebyshev Realizations of the Universal Birkhoff Theory}
\label{sec:NewBirkLeg+Cheb}

A practical realization of the universal Birkhoff theory for trajectory optimization narrows the selection of $\pi^N$ to one that satisfies the two hypotheses of \cite{newBirk-part-I}. 
In this section, we show that the Gegenbauer family of Legendre and Chebyshev grids satisfy these hypotheses, thus paving the way for a practical computational framework.  Note furthermore that the theory of \cite{newBirk-part-I} is based on an arbitrary grid and makes no assumptions of polynomial basis functions.  In other words, no claim is made here that the Legendre and Chebyshev nodes are the only grid points that meet the criteria for the Birkhoff theory to hold.

For the purposes of brevity, we rely heavily on the recent developments from \cite{newBirk-2023}, together with some well-established results in spectral methods that use the theory of orthogonal polynomials\cite{boyd,fornberg,atap,hgg-2007}. Furthermore, because this paper is a continuation of \cite{newBirk-part-I}, we use previous symbols and notation without redefining them here.

In choosing Gegenbauer nodes, we hereafter constrain the computational time interval $[\tau^a, \tau^b]$ to the $[-1, 1]$ domain (i.e., $\tau^a = -1 = - \tau^b$).  As a result, a domain transformation is needed to apply the ensuing method for a generic time interval $[t^a, t^b]$.  This detail is discussed in Section~\ref{sec:implementation}.

\subsection{Definition and $\mathcal{O}(1)$-Computation of Legendre and Chebyshev Grids}

We begin by first noting that the Chebyshev grid points can be computed at an $\mathcal{O}(1)$ computational speed using the formula, $\tau_j = -\cos(\xi_j \pi)$, where $\xi_j, \ j = 0, \ldots, N$ is determined based on the following selections\cite{hgg-2007,boyd,fornberg,trefethen-2000}:
%
\begin{subequations}\label{eq:cheb-grids}
\begin{align}
&\text{Chebyshev-Gauss-Lobatto (CGL): } & \xi_j &= \frac{ j}{N} \\
&\text{Chebyshev-Gauss-Radau (CGR): }  & \xi_j &=  \frac{2j}{2N+1}  \\
& \text{Chebyshev-Gauss (CG): }   & \xi_j &= \frac{2j +1}{2N+2}
\end{align}
\end{subequations}
%
The set of Legendre grid points are determined by the $(N+1)$-zeroes of a function $z(\tau)$, where, $z(\tau)$ and the elements of the Legendre set are selected according to the following pairings\cite{hgg-2007,boyd,fornberg,trefethen-2000}:
\begin{subequations}\label{eq:Leg-grids}
\begin{align}
&\text{Legendre-Gauss-Lobatto (LGL): } & z(\tau) &= (1- \tau^2) \dot P_N(\tau)\\
&\text{Legendre-Gauss-Radau (LGR): }  & z(\tau) &=  P_N(\tau) + P_{N+1}(\tau) \\
& \text{Legendre-Gauss (LG): }   & z(\tau) &= P_{N+1}(\tau)
\end{align}
\end{subequations}
%
where, $P_N$ is a Legendre polynomial of order $N$.  A quick examination of \eqref{eq:Leg-grids} appears to suggest that the computation of Legendre grids are significantly slower than $\mathcal{O}(1)$.  In fact, this was true up until the pioneering work of Bogaert\cite{Bogaert-2014}. Bogaert's algorithm\cite{Bogaert-2014} computes all Legendre grids at an $\mathcal{O}(1)$ computational speed.  With his advancement, all grid points denoted in \eqref{eq:cheb-grids} and \eqref{eq:Leg-grids} can be computed at an $\mathcal{O}(1)$ speed.

\subsection{A Proof of Hypotheses~1 and 2 of \cite{newBirk-part-I} for a Select Family of Gegenbauer Grids}\label{subsec:proofOfHyp1+2}
Gegenbauer polynomials (also known as ultraspherical polyonomials\cite{elbarbary}) are a subset of Jacobi polynomials that are orthogonal with respect to the weight function $(1 - \tau^2)^{\alpha-1/2}$. Setting $\alpha = 0$ yields Chebyshev polynomials (of the first kind) while $\alpha = 1/2$ generates Legendre polynomials.

To prove hypotheses~1 and 2, we need the two theorems presented below.  As noted earlier, we direct the reader to multiple sources for more complete proofs of Theorems~\ref{thm:Big2-wts} and \ref{thm:quad-error}.

\begin{theorem}[\cite{newBirk-2023}] \label{thm:Big2-wts}
\begin{enumerate}
\item[]
\item Let $\pi^N$ be any one of the family of Chebyshev grid points given by  \eqref{eq:cheb-grids}. Then the Birkhoff quadrature weights are identical to the Clenshaw-Curtis (CC) weights for CGL nodes and Fej\'{e}r's weights for CGR- and CG-grid points.
\item Let $\pi^N$ be any one of the family of Legendre grid points given by  \eqref{eq:Leg-grids}.  Then the Birkhoff quadrature weights are identical to the LGL, LGR and LG quadrature weights for the respective grid selections.
\end{enumerate}
\end{theorem}
\begin{remark}
The importance of Theorem~\ref{thm:Big2-wts} is that no new formulas are needed for the computation of the Birkhoff quadrature weights provided the universal interpolants proposed in \cite{newBirk-2023} are used.  In sharp contrast, this result does not hold for the original Birkhoff interpolants proposed by Wang et al\cite{wang}.  For instance, for an LGL grid, one needs to use the ``Gauss-Lobatto-Legendre-Birkhoff'' quadrature weights developed in \cite{GLLB}.
\end{remark}
\begin{theorem}\label{thm:quad-error}
Let $y(\cdot) \in C[-1, 1]$.  Define,
\begin{equation}
Q^N(y(\cdot)) := \sum_{i=0}^N y(\tau_i) \, w^{L/C}_i - \int_{-1}^{1} y(\tau)\, d\tau
\end{equation}
where $w^{L/C}_j, \ j = 0, \ldots, N$ are the family of Legendre/Chebyshev quadrature weights defined in Theorem~\ref{thm:Big2-wts}.  Then, given any $\epsilon > 0$, there exists an $N_\epsilon \in \mathbb{N}$ such that for all $N \ge N_\epsilon$, $\abs{Q^N(y(\cdot))} \le~\epsilon$
\end{theorem}
\begin{proof}
This theorem is well-known for the family of Legendre quadrature weights (Lobatto, Radau and Gauss).  See, for example \cite{boyd,fornberg,hgg-2007}. In the case of quadratures over a Chebyshev grid (Lobatto, Radau and Gauss), it was long thought that the results were not as accurate as those generated by Legendre node points\cite{trefethen:CC,atap}. In fact, there is widespread belief that even among the Legendre family of quadrature weights, Gauss is the best. According to Trefethen\cite{trefethen:CC}, computations performed as far back as the 1960s showed that the Clenshaw-Curtis method yielded results nearly as accurate as Gauss quadratures for the same number of node points.  More recent research\cite{trefethen:CC,fejerRef} has placed these practical observations on a firm theoretical foundation where the purported differences in quadrature errors between Legendre and Chebyshev grids (be in Lobatto, Radau or Gauss) are either nonexistent or vanishingly small.  Detailed proofs of these results can be found in \cite{atap,trefethen:CC,fejerRef,trefethen-2000}.  Note also that the statement of the current theorem is much weaker than a requirement that the Legendre/Chebyshev grid points have the same quadrature error, rather, the only requirement is that the errors be vanishingly small as $N$ increases.
\end{proof}

\begin{proposition}[Proof of Hypothesis~1 of \cite{newBirk-part-I}]\label{thm:a=b4Big2}
Let $\pi^N$ be any one of the family of Legendre/Chebyshev grids from Theorem~\ref{thm:Big2-wts}. Then, Hypothesis~1 of \cite{newBirk-part-I} holds; i.e., the $a$- and $b$-forms of the Birkhoff interpolants given in \cite{newBirk-part-I} (and \cite{newBirk-2023}) are equivalent.
\end{proposition}
\begin{proof}
Let $[-1, 1] \ni \tau \mapsto y(\tau) \in \Real$ be a given bounded differentiable function.  Then, from their defining equations given in \cite{newBirk-part-I}, it follows that,
\begin{equation}
I^N_a y(-1) = y(-1) \quad\text{and}\quad I^N_b y(1) = y(1)
\end{equation}
Hence,
\begin{equation}\label{eq:Ib-Ia=byDef}
I^N_b y(\tau) - I^N_a y(\tau) = y(1) - y(-1) + \sum_{j=0}^N \dot y(\tau_j) \left(B_j^b(\tau) - B_j^a(\tau)\right)
\end{equation}
From Proposition~3 of \cite{newBirk-2023} (see also Lemma~2 in \cite{newBirk-part-I}), \eqref{eq:Ib-Ia=byDef} can be written as,
\begin{equation}\label{eq:Ib-Ia=proof-1}
I^N_b y(\tau) - I^N_a y(\tau) = y(1) - y(-1) - \sum_{j=0}^N \dot y(\tau_j) w^B_j
\end{equation}
%
From Theorem~\ref{thm:Big2-wts}, \eqref{eq:Ib-Ia=proof-1} can be written as,
\begin{equation}\label{eq:Ib-Ia=proof-2}
I^N_b y(\tau) - I^N_a y(\tau) = y(1) - y(-1) - \sum_{j=0}^N \dot y(\tau_j) w^{L/C}_j
\end{equation}
From Theorem~\ref{thm:quad-error}, \eqref{eq:Ib-Ia=proof-2} can be written as,
\begin{equation}\label{eq:Ib-Ia=proof-3}
I^N_b y(\tau) - I^N_a y(\tau) = y(1) - y(-1) - \int_{-1}^1 \dot y(\tau)\, d\tau + \epsilon^N
\end{equation}
where $\epsilon^N$ is vanishingly small as $N \to \infty$. The result then follows by evaluating the integral in \eqref{eq:Ib-Ia=proof-3}.
\end{proof}
\begin{proposition}[Proof of Hypothesis~2 of \cite{newBirk-part-I}]
Theorem~\ref{thm:quad-error} is a proof of Hypothesis~2 of \cite{newBirk-part-I} for the choices of $\pi^N$ denoted in Theorem~\ref{thm:Big2-wts}.
\end{proposition}

\subsection{Impact of Selecting Legendre/Chebyshev- Gauss, Gauss-Radau and Gauss-Lobatto Grids}
The results of subsection~\ref{subsec:proofOfHyp1+2} suggest that there is virtually no difference in selecting any one of the six families of grid points denoted in Theorem~\ref{thm:Big2-wts}.  Although this is true in principle, we note the following caveats that are specific to computational optimal control:
\begin{enumerate}
\item \textbf{Legendre/Chebyhsev-Gauss grid:} A Gauss grid does not include the boundary points, $-1$ and $1$.  Consequently, if the boundary conditions are imposed at the boundary points, the controls at the boundaries are unavailable. Hence, to perform a feasibility test (see Section~\ref{sec:implementation}), some form of control extrapolation is needed.
\item \textbf{Legendre/Chebyhsev-Gauss-Radau grid:} A Radau grid does not include one of the boundary points.  Thus, either the initial- or final-value of the control is unavailable.  Extrapolation is needed in both cases.
\item \textbf{Legendre/Chebyhsev-Gauss-Lobatto grid:} A Lobatto grid includes both boundary points; hence, both the initial and final values of the controls are available without introducing any adhoc extrapolation rules.
\end{enumerate}
Although the caveats noted above are relatively minor with regards to a Birkhoff implementation, note that these small differences are amplified to a detrimental effect with regards to Gauss or Radau selections when a Lagrange PS method is chosen for trajectory optimization\cite{advances,gongProbTN,PSReview-ARC-2012}. This is the main reason why Gauss/Radau grids have not been selected for flight implementations\cite{PSReview-ARC-2012,flight:TRACE-2014,zpm:IEEE,TRACE-IEEE-Spectrum}.

\begin{remark}\label{rem:classification}
Following \cite{wang,newBirk-2023,newBirk-part-I}, we classify PS methods based on the the type of interpolants used; i.e., Lagrange or Birkhoff. Secondary classification follows the well-established convention\cite{boyd,fornberg,trefethen-2000} of basing it on the type of computational basis functions used; e.g., Legendre, Chebyshev etc.  A tertiary level of classification is based on the type of grid used with the associated computational basis function; e.g., Lobatto, Radau, Gauss.
\end{remark}




\section{Legendre/Chebyshev Implementations of Birkhoff Optimal Control Theory}\label{sec:implementation}

A Legendre/Chebyshev implementation of the universal Birkhoff theory for optimal control\cite{newBirk-part-I} can now be performed by incorporating the results of Section~\ref{sec:NewBirkLeg+Cheb} with those from \cite{newBirk-part-I} and \cite{newBirk-2023}.  Additional details of implementation follow from relevant optimal control programming techniques that are extensively described elsewhere\cite{PSReview-ARC-2012,spec-alg,DIDO:arXiv,fastmesh,furtherResults,auto-knots}.  Hence, we limit the scope of this section to just those details that are specific only to the new Birkhoff method.

\subsection{A Generic Optimal Control Problem Formulation}
A fairly generic optimal control problem may be described in the form of following compact format\cite{ross-book}:
%
\begin{eqnarray}
&\bx(t) \in \real{N_x}, \quad \bu(t) \in \real{N_u}, \quad t \in [t^a, t^b], \quad \bp \in \real{N_p}   & \nonumber\\
& (\textsf{$\bP$}) \left\{
\begin{array}{lrl}
\emph{Minimize } & J[\bx(\cdot), \bu(\cdot), t^a, t^b, \bp] =& E(\bx(t^a),\bx(t^b), t^a, t^b, \bp)\\
& & \quad \displaystyle + \int_{t^a}^{t^b} F(\bx(t), \bu(t), t, \bp)\, dt \\
\emph{Subject to}& \dot\bx(t) =& \bff(\bx(t), \bu(t), t, \bp)  \\
& \be^L \le & \be(\bx(t^a), \bx(t^b), t^a, t^b, \bp) \le \be^U \\
& \bh^L \le & \bh(\bx(t), \bu(t), t, \bp) \le \bh^U
\end{array} \right. & \label{eq:probPbold}
\end{eqnarray}
%
where, the symbol $N_{(\cdot)} \in \mathbb{N}$ represents the number of variables or equations implied by its subscript.  We assume $t^b > t^a$ and $\bp$ is a static optimization parameter that must be determined jointly with $\bxf$ and $\buf$ in addition to the clock times, $t^a$ and $t^b$.  The cost functional $J$ is given in the ``Bolza'' or standard form\cite{ross-book} where, $E$ and $F$ are the the endpoint- and running-cost functions respectively. The running cost, $F$, pairs with the dynamics data function, $\bff$, via the adjoint covector (costate), $\blam$, to form the (Pontryagin) Hamiltonian function,
\begin{equation}
H(\blam, \bx, \bu, t, \bp) := F(\bx, \bu, t, \bp) + \blam^T\bff(\bx, \bu, t, \bp)
\end{equation}
The function, $E$, pairs with the endpoint constraint function, $\be$, via the events covector, $\bnu$, to form the endpoint Lagrangian\cite{ross-book}:
\begin{equation}
\overline{E}(\bnu, \bx^a, \bx^b, t^a, t^b, \bp):= E(\bx^a, \bx^b, t^a, t^b, \bp) + \bnu^T \be(\bx^a, \bx^b, t^a, t^b, \bp)
\end{equation}
Likewise, the path constraint function, $\bh$, pairs with $H$ via the path covector, $\bmu$, to form the Lagrangian of the Hamiltonian; i.e., the Lagrangian associated with the Hamiltonian minimization condition (HMC)\cite{ross-book},
\begin{equation}
\overline{H}(\bmu, \blam, \bx, \bu, t, \bp) := H(\blam, \bx, \bu, t, \bp)  + \bmu^T  \bh(\bx, \bu, t, \bp)
\end{equation}
The superscripts $L$ and $U$ in \eqref{eq:probPbold} denote their lower and upper bounds on the values of the functions $\be$ and $\bh$.

Additional cost functionals and constraints may be added to further generalize \eqref{eq:probPbold}; however, it will be apparent later that such additions can be easily incorporated into the Birkhoff theory using the fundamentals covered in \cite{newBirk-part-I} and this section.  For similar reasons, incorporating ``phases'' in generalizing \eqref{eq:probPbold} follows in quite a straightforward manner based on the results of \cite{acc:hybrid,knots,spec-alg,auto-knots,hybrid:jgcd}.  As noted in \cite{newBirk-part-I} and elsewhere\cite{fastmesh,furtherResults,advances}, such details, while important, are more distracting to the presentation of the core fundamentals.

\subsection{Legendre/Chebyshev Details for Birkhoff Optimal Control Programming}\label{subsec:BirkProgDetails}
To apply the theory from \cite{newBirk-part-I} for a selection of Legendre/Chebyshev grid points to solve Problem~$(\bP)$, two additional ingredients are necessary: (1) an invertible and differentiable domain transformational mapping (i.e., a diffeomorphism)  from the $[-1, 1]$ to $[t^a, t^b]$, and (2) explicit computational formulas for producing $\bB^\theta(\pi^N), \ \theta \in \set{a, b}$  and $\bw_B(\pi^N)$ for the choices of $\pi^N$ denoted in Theorem~\ref{thm:Big2-wts}.

\subsubsection{Domain Transformation}
Let, $\Gamma$ be a differentiable, invertible function with arguments $ (\tau, t^a, t^b; \bp_t)$, such that\cite{PSReview-ARC-2012,ross-book,DIDO:arXiv},
\begin{equation}\label{eq:t=Gamma}
t = \Gamma(\tau, t^a, t^b; \bp_t)
\end{equation}
where,  $t^a = \Gamma(-1, t^a, t^b; \bp_t)$,  $t^b = \Gamma(1, t^a, t^b; \bp_t)$ and $\bp_t \in \real{N_{p_t}}$ is a parameter that can be selected a priori or optimized.  If $\bp_t$ is optimized, the subsequent grid is adaptive\cite{PSReview-ARC-2012,kt-map}; if $\Gamma$ is nonlinear, the grid can be designed to support anti-aliasing\cite{bellman-low-t,bellman-conf} and assist in the resolution of the system trajectory\cite{ross-book,kt-map,bellman-conf}.  The simplest choice of $\Gamma$ for a finite horizon problem is the affine function,
\begin{equation}\label{eq:linear-tX}
t =\Gamma(\tau, t^a, t^b; \bp_t) \equiv \Gamma_{aff}(\tau, t^a, t^b) := \left(\frac{t^b - t^a}{2} \right)\tau  +  \left(\frac{t^b + t^a}{2} \right)
\end{equation}
See \cite{ross-book,LGR-rtoc,Radau-JGCD,Radau-GNC05} on how to select $\Gamma$ for infinite horizon problems.  With a suitable choice of $\Gamma$, Problem~$(\bP)$ is easily transformed to the computational domain $[-1, 1]$ by substituting \eqref{eq:t=Gamma} in \eqref{eq:probPbold}.  The integral and the differential equations in \eqref{eq:probPbold} transform according to,
\begin{subequations}\label{eq:t2tau-trans}
\begin{align}
\int_{t^a}^{t^b} F(\bx(t), \bu(t), t, \bp)\, dt &=\int_{-1}^{1} \gamma\, F(\bx(\tau), \bu(\tau), \Gamma, \bp)\, d\tau  \\
\frac{d\bx(\tau)}{d\tau} &=\gamma\, \bff(\bx(\tau), \bu(\tau), \Gamma, \bp) \label{eq:xdot=f-transformed}
\end{align}
\end{subequations}
where, $\gamma =  \left(\frac{d\Gamma}{d\tau}\right)$.  In \eqref{eq:t2tau-trans} we have abused notation for expediency.  The symbols $\bx(\tau)$ and $\bu(\tau)$ should be written more technically as $\bx\big(\Gamma(\tau, t^a, t^b; \bp_t)\big)$ and $\bu\big(\Gamma(\tau, t^a, t^b; \bp_t)\big)$ respectively.  Likewise, the dependence of $\gamma$ on $t^a, t^b$ and possibly $\bp_t$ is suppressed for convenience.

\subsubsection{Computation of $\bB^\theta(\pi^N)$ and $\bw_B(\pi^N)$}
Efficient computational procedures for computing $\bB^\theta(\pi^N), \ \theta \in \set{a, b}$ over a generic grid, $\pi^N$, are provided in \cite{newBirk-2023}.  In particular, \cite{newBirk-2023} provides explicit computational formulas for $\bB^\theta(\pi^N)$ when $\pi^N$ is chosen in accordance with the selection of grids denoted in Theorem~\ref{thm:Big2-wts}.

To compute $\bw_B(\pi^N)$, the following options are available:
\begin{enumerate}
\item If $\pi^N$ is chosen to be an LGL or a CGL grid, then the last row of $\bB^a$ turns out to be exactly equal to $\bw_B$\cite{newBirk-2023}; hence, no new computations are necessary.
\item If $\pi^N$ is chosen to be a non-Lobatto grid (i.e, Radau or Gauss), then additional computations must be performed to produce $\bw_B$.  The software, \texttt{chebfun}\cite{chebfun-2014}, implements fast and efficient formulas for computing these weights (cf.~Theorem~\ref{thm:Big2-wts}).
\end{enumerate}
To verify and validate these computation, the following checks may be performed:
\begin{enumerate}
\item $B^a_j(\tau_k) - B^b_j(\tau_k) = w^B_j, \quad j, k = 0, \ldots, N$; see Proposition~3 in \cite{newBirk-2023}.
\item If $\pi^N$ is selected to be an LGL or a CGL grid, then, $\bB^a = -I^\diamondsuit \bB^b I^\diamondsuit $, where, $I^\diamondsuit$ is the exchange matrix; see Proposition~5 in \cite{newBirk-2023}.
\end{enumerate}

\subsubsection{Matrix-Vector Mathematical Programming Problem Formulation}
As noted elsewhere\cite{DIDO:arXiv,furtherResults,ross:shooting-arXiv,bilevel-2014,bilevel-2016,Hamil-mesh-2016}, discretizing Problem~$(\bP)$ (by any valid method) generates an information-rich mathematical programming problem that is better described as a Hamiltonian programming\cite{DIDO:arXiv,Hamil-mesh-2016,ross:shooting-arXiv} problem rather than a generic nonlinear programming problem. To preserve this Hamiltonian structure, we follow \cite{DIDO:arXiv,fastmesh} and define three  matrices according to:
\begin{subequations}
\begin{align}
\bX &:= [\bx_0, \bx_1, \ldots, \bx_N] \quad \in \real{N_x \times (N+1)}  \label{eq:xmat}\\
\bV &:= [\bv_0, \bv_1, \ldots, \bv_N] \quad \in \real{N_x \times (N+1)} \\
\bU &:= [\bu_0, \bu_1, \ldots, \bu_N] \quad \in \real{N_u \times (N+1)}  \label{eq:umat}
\end{align}
\end{subequations}
where, $\bx_k \in \real{N_x}, \bv_k \in \real{N_x}, \bu_k \in \real{N_u}, \ k = 0, 1, \ldots, N$. This format 
allows us to preserve the structural nature of \eqref{eq:probPbold}; for instance, the path constraint function can be evaluated over the transformed grid, $\Gamma^N = \Gamma(\pi^N, t^a, t^b, \bp)$, by overloading $\bh$ according to\cite{ross-book,DIDO:arXiv},
\begin{equation}
\bh(\bx(t), \bu(t), t, \bp) \mapsto \bh(\bX, \bU, \Gamma^N, \bp) \quad \in \real{N_h \times (N+1)}
\end{equation}
Using similar constructs and the framework of \cite{newBirk-part-I}, it is straightforward to show that the $a$-form of the universal Birkhoff discretization of Problem $(\bP)$ can be framed as:
\begin{eqnarray}
&\bX \in \real{N_x \times (N+1)}, \quad \bU \in \real{N_u \times (N+1)}, \quad \bV \in \real{N_x \times (N+1)}\nonumber \\
& t^a \in \Real, \quad t^b \in \Real, \quad \bx^a \in \real{N_x}, \quad \bx^b \in \real{N_x}  & \nonumber\\
& \left(\textsf{$\bP_a(\pi^N)$}\right) \left\{
\begin{array}{lrl}
\emph{Minimize } & J^N[\bX, \bU, \bV, \bx^a, \bx^b, t^a, t^b] :=& E(\bx^a, \bx^b, t^a, t^b)\\
& & + \ \gamma\, \bw_B^T(\pi^N) \, F(\bX, \bU, \Gamma^N, \bp)  \\[1em]
\emph{Subject to} &\bA_a(\pi^N) \left[
                                \begin{array}{c}
                                  \bX^T \\
                                  \bV^T \\
                                \end{array}
                              \right]
  -& \bC^N_a \left[
                                \begin{array}{c}
                                  \bx^a \\
                                  \bx^b \\
                                \end{array}
                              \right]  = \bzero\\[2em]
& \bV =& \gamma\, \bff(\bX, \bU, \Gamma^N, \bp)  \\
& \be^L \le \be(\bx^a, \bx^b, t^a, t^b, \bp)  \le & \be^U \\
& \bh^L \le \bh(\bX, \bU, \Gamma^N, \bp)  \le & \bh^U
\end{array} \right. & \label{eq:ProbbPNa}
\end{eqnarray}
where, $\bA_a(\pi^N)$ is an $\big((N+2) \times 2(N+1)\big)$-matrix that depends upon the choice of $\pi^N$,
\begin{align}
\bA_a(\pi^N) := \left[
                     \begin{array}{cc}
                      I^N & -\bB^a(\pi^N) \\
                    \bzero^T & \bw_B^T(\pi^N)  \\
                                \end{array}
                              \right]
\end{align}
$\bzero$ is an $\big((N+1)\times 1\big)$ vector of zeros, and $\bC^N_a$ is an $\big((N+2) \times 2\big)$-matrix given by,
\begin{align}
\bC^N_a := \left[
                     \begin{array}{cc}
                      \bb & \bzero \\
                    -1 & 1  \\
                                \end{array}
                              \right]
\end{align}
The symbol $\bP_a(\pi^N)$ used in \eqref{eq:ProbbPNa} is the vectorized version of $P^N_a$ defined in \cite{newBirk-part-I}.  In \eqref{eq:ProbbPNa} it is written more elaborately as $\bP_a(\pi^N)$ to signify its dependence on the choice of $\pi^N$.  See also \cite{furtherResults} for additional details on efficient matrix methods for organizing variables.
\begin{remark}
Similar to the development of \eqref{eq:ProbbPNa}, other forms of $\bP_\theta(\pi^N), \ \theta \in \set{b, a^*, b^*}$\cite{newBirk-part-I} can be generated.  Because the process is quite similar, they are not detailed here for brevity but depicted in Fig.~\ref{fig:BirkhoffPStypes}.  It is apparent from Fig.~\ref{fig:BirkhoffPStypes} that a multiplicity of possibilities can be easily generated using the results of Section~\ref{sec:NewBirkLeg+Cheb} and the theory developed in \cite{newBirk-part-I}.
%
\begin{figure}[h!]
      \centering
      {\parbox{0.9\textwidth}{
      \centering
      {\includegraphics[width = 0.7\textwidth]{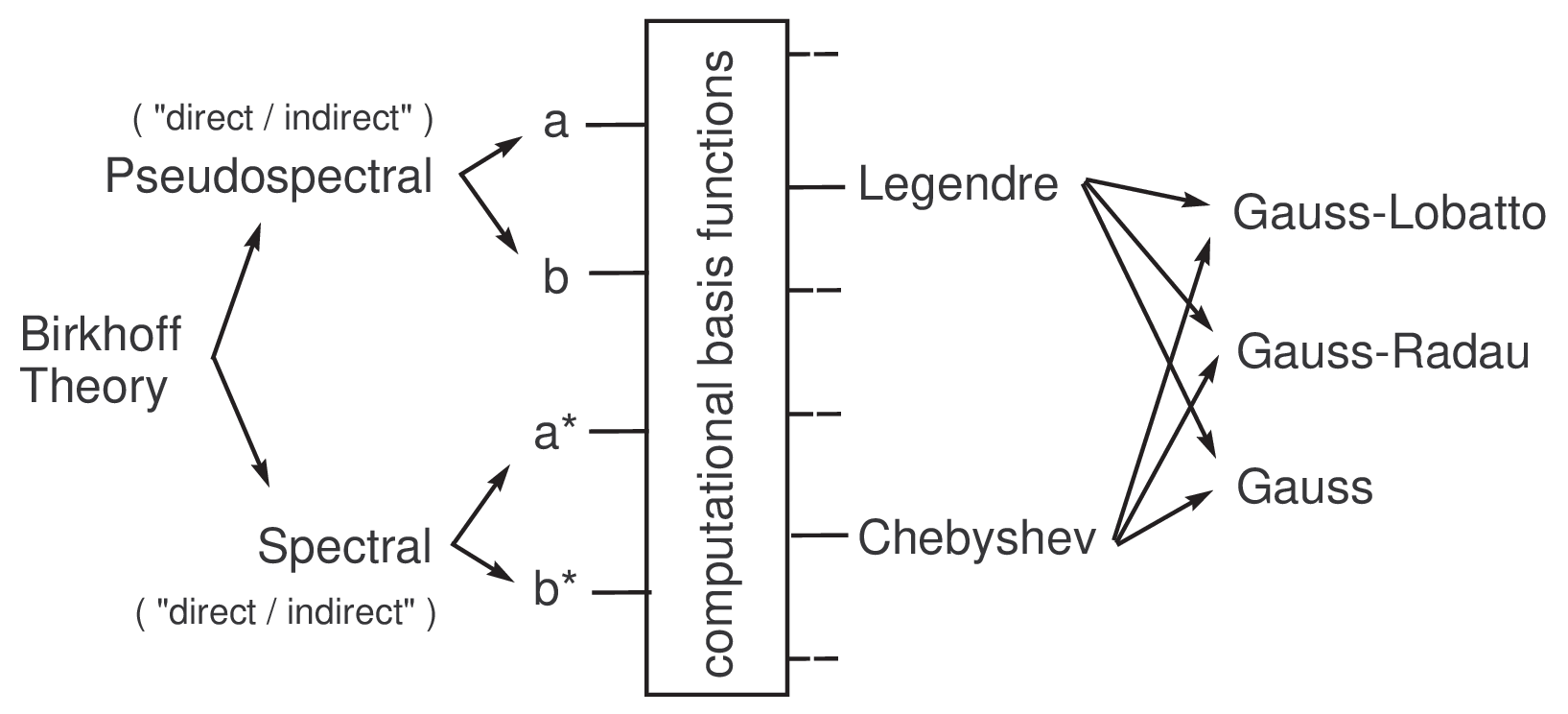}}
      \caption{\textsf{Schematic for the flow of the new Birkhoff theory for trajectory optimization: from its origin to various primal/dual implementations; figure adapted from \cite{DIDO:arXiv}.}}\label{fig:BirkhoffPStypes}
      }
      }
\end{figure}
%
See also Remark~\ref{rem:classification} in the context of Fig.~\ref{fig:BirkhoffPStypes}.
\end{remark}

Once a computable version of Problem~$(\bP_\theta(\pi^N)), \ \theta \in \set{a, b, a^*, b^*}$ is formulated, then the entire power of the universal Birkhoff theory\cite{newBirk-part-I} can be brought to bear on solving Problem~$(\bP)$.  A candidate solution obtained by such a process can then be subjected to a rigorous battery of Birkhoff-independent verification and validation tests.  These tests are described in the next section and illustrated via sample problems and solutions in Section~\ref{sec:examples}.

\section{Birkhoff-Independent Verification and Validation Tests}\label{sec:IV&V}
Suppose one is presented with a candidate optimal control trajectory, $t \mapsto \bu(t) \in \real{N_u}$, that purports to solve Problem~$(\bP)$.  Regardless of how $\bu(\cdot)$ was obtained, it is obvious that an optimal trajectory must satisfy the necessary conditions for optimality (otherwise, these conditions would be unnecessary).  A precursor to optimality is primal feasibility.  This is simply because if the candidate control trajectory is not feasible, it is certainly not optimal.

\subsection{Primal Feasibility Tests}\label{subsec:primal-feasibility}

In any grid-based computation, the control signal is given by the set of points,
\begin{equation}\label{eq:u=grid-soln}
\set{\big(t_i, \bu(t_i)\big) \in \Real \times \real{N_u}, \ i = 0, \ldots N }
\end{equation}
Frequently, but not always, this is associated with a state solution over the same grid:
\begin{equation}\label{eq:x=grid-soln}
\set{\big(t_i, \bx(t_i)\big) \in \Real \times \real{N_x}, \ i = 0, \ldots N }
\end{equation}
A claim of primal feasibility based solely on the pairs $\big(\bx(t_i), \bu(t_i), \ i = 0, \ldots, N\big)$ satisfying all the the constraints associated with Problem~$(\bP)$ (cf.~\eqref{eq:probPbold}, with $\bp$ temporarily ignored for convenience) is frequently false\cite{ross-book,bellman-low-t,bellman-conf}.  To clarify this statement, we use the following definition:
%
\begin{definition}\label{def:primal-feasibility}
A state-control function pair, $t \mapsto \big( \bx^\sharp(t), \bu^\sharp(t) \big) \in \real{N_x} \times \real{N_u}$ and a parameter $\bp^\sharp$ are said to be primal feasible if they satisfy all the constraints associated with Problem~$(\bP)$.
\end{definition}
%
Note that Definition~\ref{def:primal-feasibility} is not based on \eqref{eq:u=grid-soln} and \eqref{eq:x=grid-soln}. To obtain a \emph{certificate of feasibility} for \eqref{eq:u=grid-soln}, one must first generate a continuous-time control signal $t \mapsto \bu^\sharp(t)$.  This can be easily done by using the formula,
\begin{equation}\label{eq:uSharp=cont}
\bu^\sharp(t) = \sum_{i=0}^N \bu(t_i)\zeta_i(t)
\end{equation}
where, $\zeta(t)$ is any valid interpolating function; i.e., one that satisfies the condition $\zeta_i(t_j) = \delta_{ij}$.  In practice, it is highly desirable for $\zeta(t)$ to be agnostic to the method on how \eqref{eq:u=grid-soln} was created so that there is an additional layer of independence in the tests to be performed. See, for example \cite{TEI-JGCD-2011} for a detailed discussion on choosing $\zeta_i(t)$ for a practical problem pertaining to a trans-Earth Lunar mission.  A state trajectory, $t \mapsto \bx^\sharp(t)$, associated with \eqref{eq:uSharp=cont} can now be obtained quite easily by solving the initial value problem (IVP),
\begin{equation}\label{eq:IVP-V&V}
\dot\bx^\sharp = \bff(\bx^\sharp, \bu^\sharp(t)) \equiv \bg(\bx^\sharp, t), \quad \bx^\sharp(t^a) = \bx^a
\end{equation}
Thus, a procedure for a primal feasibility test is as follows:
\begin{enumerate}
\item Solve Problem~$(\bP_\theta(\pi^N))$ to generate \eqref{eq:u=grid-soln}.
\item Select a $\zeta(t)$ (typically, piecewise linear) to construct \eqref{eq:uSharp=cont}.
\item Solve the IVP given by \eqref{eq:IVP-V&V} to generate $t \mapsto \bx^\sharp(t)$.
\item If $t \mapsto \big( \bx^\sharp(t), \bu^\sharp(t) \big)$ is primal feasible (cf.~Definition~\ref{def:primal-feasibility}), then we declare the control solution given by \eqref{eq:u=grid-soln} to be feasible.
\end{enumerate}
%
\begin{remark}\label{rem:V&V-primal}
For flight implementation, it is irrelevant how \eqref{eq:u=grid-soln} was generated. In typical flight operations\cite{zpm:IEEE,bhatt:opm,TRACE-IEEE-Spectrum,flight:TRACE-2014,Minelli-AAS,PSReview-ARC-2012},
\eqref{eq:uSharp=cont} is used with some $\zeta_i(t)$ (e.g.~piecewise linear or zero-order-hold) to solve \eqref{eq:IVP-V&V} (via, for example, a fourth-order Runge-Kutta integrator) to perform an independent verification of (primal) feasibility. Such tests are performed in all examples reported in Section~\ref{sec:examples}.
\end{remark}
\begin{remark}
It is worth emphasizing at this point that a grid-based control solution (in the format of \eqref{eq:u=grid-soln}) obtained by a Birkhoff method is not necessarily of polynomial origin even if the underlying computation of a Birkhoff matrix is based on orthogonal polynomials.  In fact, as explained in Section~III of \cite{newBirk-part-I}, $\bu(t_i)$ is frequently not of polynomial origin.
\end{remark}
%
See \cite{ross-book} and \cite{DIDO:arXiv} for additional details and the impact of using \eqref{eq:uSharp=cont} and \eqref{eq:IVP-V&V} with regards to practically achievable tolerances on constraint violations.

\subsection{Optimality Tests via Pontryagin's Principle}\label{subsec:dual-feasibility}
A necessary condition is based on the principle that a candidate primal-feasible solution has been found. Suppose that the tests indicated in subsection~\ref{subsec:primal-feasibility} have been successfully completed. If this feasible solution is purportedly optimal, then it must necessarily satisfy a set of conditions for optimality. This set of necessary conditions is given by Pontryagin's principle\cite{ross-book}.  The necessary conditions for Problem~$(\bP)$ are given by\cite{ross-book}:
%
\begin{subequations}\label{eq:HAMVET}
\begin{align}
&\text{Hamiltonian Minimization Condition: } && \bu(t) = \arg\min_{\bu(t) \in \U(t, \bx(t))} H(\blam(t), \bx(t), \bu(t), t, \bp) \label{eq:HMC}\\
&\text{Adjoint equation: } && -\dot\blam(t) = \partial_{\bx} \overline{H}(\bmu(t), \blam(t), \bx(t), \bu(t), t, \bp)\label{eq:adjoint}\\
&\text{Hamiltonian Value Conditions: } && \mathcal{H}[@ t^\theta] = \pm \partial_{t^\theta} \overline{E}(\bnu, \bx^a, \bx^b, t^a, t^b, \bp)\quad \theta \in \set{a, b} \\
&\text{Hamiltonian Evolution Equation: } && \frac{d \mathcal{H}(\blam(t), \bx(t), t, \bp)}{dt}  \nonumber\\
& && \qquad = \frac{\partial \overline{H}(\bmu(t), \blam(t), \bx(t), \bu(t), t, \bp)}{\partial t}\\
&\text{Transversality Condition: }  && \blam(t^\theta) = \mp \partial_{\bx^\theta} \overline{E}(\bnu, \bx^a, \bx^b, t^a, t^b, \bp)\quad \theta \in \set{a, b} \label{eq:TVC}\\
&\text{Optimal Parameter Condition: } && \partial_{\bp}\overline{E}(\bnu, \bx^a, \bx^b, t^a, t^b, \bp) \nonumber\\
& && \quad + \int_{t^a}^{t^b} \partial_{\bp} \overline{H}(\bmu(t), \blam(t), \bx(t), \bu(t), t, \bp)\, dt = \bzero
\end{align}
\end{subequations}
%
where, $\bnu$ and $\bmu(t)$ satisfy the complementarity conditions denoted symbolically as\cite{ross-book},
\begin{align}\label{eq:complementarity}
\bnu \dagger \be(\bx^a, \bx^b, t^a, t^b, \bp), && \bmu(t) \dagger \bh(\bx(t), \bu(t), t, \bp) \quad\forall\ t \in [t^a, t^b]
\end{align}

As explained in Section~\ref{sec:Intro}, it is necessary for certain operators to be commutative because an absence of commutativity may generate the converged but wrong answer; see \cite{ross:roadmap-2005} for details. The required commutativity between dualization and discretization is indeed the CMP and is illustrated in Fig.~3 of \cite{newBirk-part-I}. Because the Legendre/Chebyshev grids denoted in Theorem~\ref{thm:Big2-wts} satisfy the two hypotheses required for the universal Birkhoff theory to hold, the covector mapping theorems of \cite{newBirk-part-I} hold for this selection of node points. Consequently, all the conditions stipulated in \eqref{eq:HAMVET} and \eqref{eq:complementarity} must hold if a candidate control solution given by \eqref{eq:u=grid-soln} is to be optimal. Such primal-feasibility- and optimality tests are illustrated in Section~\ref{sec:examples}.


\section{Illustrative Problems and Solutions}\label{sec:examples}
As implied by the theoretical foundations presented in \cite{newBirk-part-I}, a computational method based on the universal Birkhoff theory should have no major technical roadblocks in solving a wide variety of optimal control problems. The main caveat for the theory to hold is the two hypotheses put forth in \cite{newBirk-part-I}.  Having shown in Section~\ref{sec:NewBirkLeg+Cheb} that these hypotheses are indeed satisfied for the family of Legendre and Chebyshev grids denoted in Theorem~\ref{thm:Big2-wts}, we now demonstrate various aspects of the numerical performance of a Birkhoff PS method via the implementations described in Section~\ref{sec:implementation}.  For the purpose of brevity, we limit the presentation of the numerical performance to CGL grids while noting that similar results were obtained for the LGL grid. 

In the following subsections, we solve a selection of problems to demonstrate the following features that are important in solving application problems in aerospace and mechanical engineering:
\begin{enumerate}
\item Production of a discontinuous optimal control solution without the need for introducing PS knots\cite{knots,spec-alg,auto-knots} for mesh refinement; 
\item Solution to a state-constrained problem that involves shifted Dirac-delta-type functions (of atomic measure) for its path covector trajectory; and,
\item Generation of a highly oscillatory trajectory, with $N > 1000$,
without facing technical roadblocks involving high condition numbers, round-off errors or instability issues.
\end{enumerate}

\subsection{A Trajectory Optimization Problem With Control Jumps}
\begin{figure}[h!]
      \centering
\begin{subfigmatrix}{2} 
\subfigure[Birkhoff-theoretic bang-bang thrust control profile and its corresponding switching covector trajectory, $t \mapsto \mu(t)$.]{\includegraphics{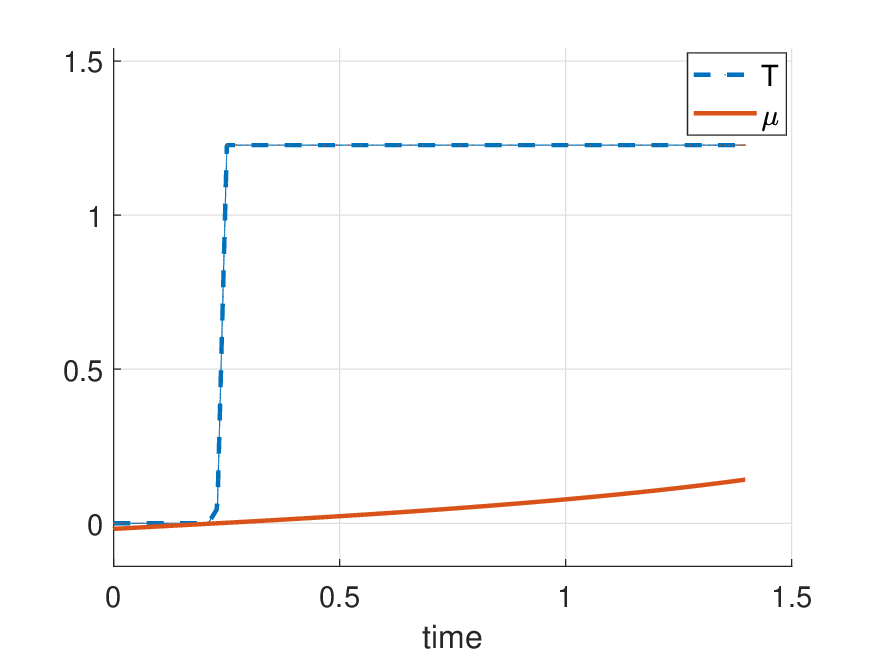}}
\subfigure[Primal feasible solution obtained via solving \eqref{eq:IVP-V&V} overlaid with a 20-node solution.]{\includegraphics{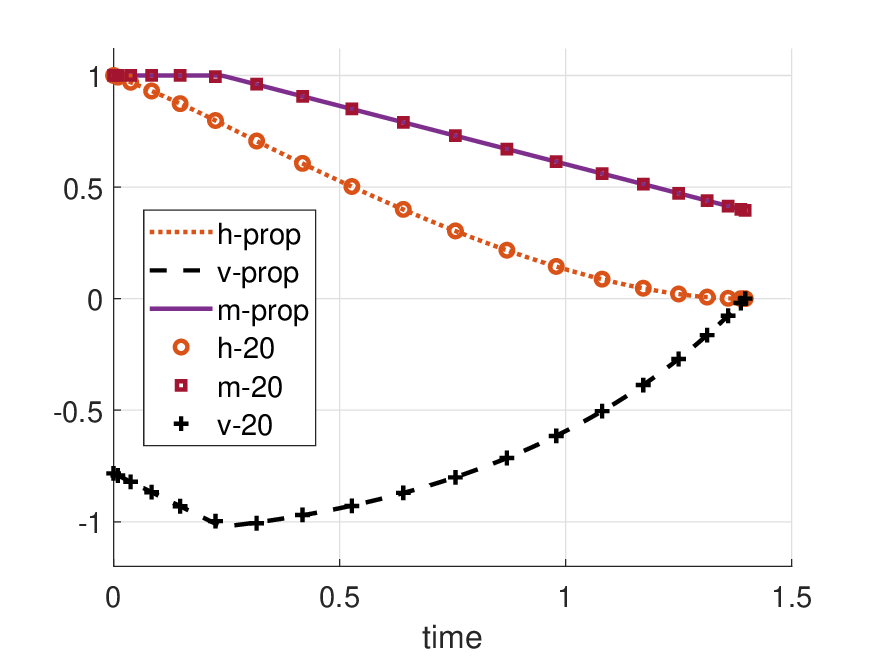}}
\end{subfigmatrix}
\caption{\textsf{Candidate optimal solution to Problem~$(ML_1)$.}}\label{fig:ML1-primals}
\end{figure}

The collection of moon-landing problems described in \cite{ross-book} are nonlinear trajectory optimization problems with control discontinuities.  A particular one-degree-of-freedom version of this problem\cite{meditch}, with moon's gravity normalized to unity, is posed in \cite{ross-book} as:
\begin{eqnarray}
& \bx := (h, v, m) \in \real 3, \quad \bu := T \in \Real   & \nonumber \\ [1em]
&(ML_{1})\left\{
\begin{array}{lrl}
\textsf{Minimize } & J[\bx(\cdot), \bu(\cdot), t_f] :=& \displaystyle \int_{t^a}^{t_f} \frac{T(t)}{2.349}\, dt
\\
\textsf{Subject to}
& \dot h(t) =& v(t) \\
& \dot v(t) =& - 1 + \frac{T(t)}{m(t)} \\
&\dot m(t) =& - \frac{T(t)}{2.349}\\
&t_0 = & 0 \\
&\big(h(t_0), v(t_0), m(t_0)\big) =& \big(1, -0.783, 1\big)\\
&(h_f, v_f) =& (0, 0)\\
& 0 \le  T \le &  1.227
\end{array} \right. &\label{eq:Prob-ML1}
\end{eqnarray}
A set of near-universal checkable conditions for optimality is the KKT conditions associated with minimizing the Hamiltonian given by \eqref{eq:HMC}. For Problem~$(ML_1)$, the complementarity condition (see \eqref{eq:complementarity}) associated with minimizing the Hamiltonian is given by\cite{ross-book},
\begin{equation}\label{eq:ML-V&V2}
\mu(t) \left\{
         \begin{array}{ll}
           \le 0, & \hbox{if } T(t) = 0 \\
           \ge 0, & \hbox{if } T(t) = 1.227
         \end{array}
       \right.
\end{equation}
where $\mu(t) \in \Real $ is the instantaneous value of the path covector associated with the thrust -limiting path constraint, $0 \le T(t) \le 1.227$.  It can be shown\cite{ross-book} that the optimal solution to this problem comprises only bang-bang controls. A linearly-interpolated 80-node Birkhoff solution for the thrust profile is shown in Fig.~\ref{fig:ML1-primals}.a) along with a similar 80-node plot of the path covector trajectory, $t \mapsto \mu(t)$.
%
%
From this plot, it is apparent that:
\begin{enumerate}
\item The on-off nature of the optimal thrust solution is well-captured by a Birkhoff solution (without introducing traditional mesh refinement techniques\cite{knots,spec-alg,auto-knots};
\item There is no Gibbs phenomenon in the bang-bang nature of the thrust profile; and
\item The optimality condition given by \eqref{eq:ML-V&V2} is satisfied.

\end{enumerate}
Fig.~\ref{fig:ML1-primals}.b) shows the state trajectories obtained by integrating (using \texttt{ode45} in MATLAB) the differential equations using the Birkhoff thrust profile (of Fig.~\ref{fig:ML1-primals}.a)).  These are denoted as $(h, v, m)$-prop in Fig.~\ref{fig:ML1-primals}.b).  These plots are in accordance with solving \eqref{eq:IVP-V&V} for an ``independent'' verification of primal feasibility; see Remark~\ref{rem:V&V-primal}. Also shown in Fig.~\ref{fig:ML1-primals}.b) are the Birkhoff-generated state solutions for $N = 20$ denoted as $(h, v, m)$-20. It is apparent from this figure that the Birkhoff solution for as little as $N = 20$ has indeed converged, thus demonstrating its fast convergence property.

\begin{figure}[h!]
      \centering
      {\parbox{0.9\textwidth}{
      \centering
      {\includegraphics[width = 0.6\textwidth]{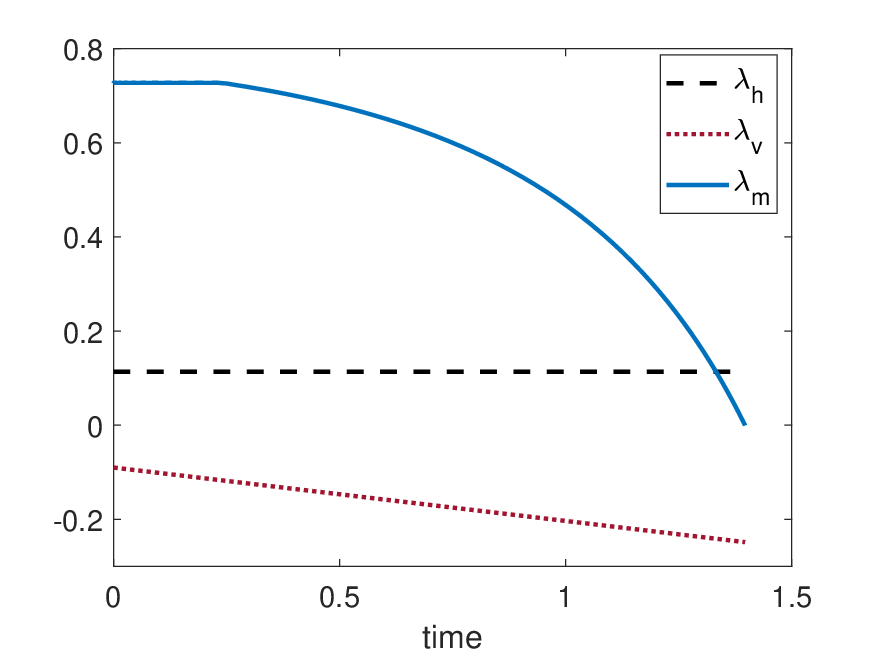}}
      \caption{\textsf{Birkhoff-theory-generated costates for the moon-landing problem given by \eqref{eq:Prob-ML1} }}\label{fig:costates-ML1}
      }
      }
\end{figure}

Additional checkable conditions can be derived for Problem~$(ML_1)$ by an application of the necessary conditions given in \eqref{eq:HAMVET}.  Using \eqref{eq:adjoint} and \eqref{eq:TVC}, it can be shown\cite{ross-book} that the costates satisfy the (checkable) conditions,
\begin{subequations}\label{eq:ML-V&V1}
\begin{align}
\lambda_h(t) &= c_1 \\
\lambda_v(t) & = - c_1 t + c_2 \\
\lambda_m(t_f) & = 0
\end{align}
\end{subequations}
where $c_1$ and $c_2$ are some constants. A plot of the Birkhoff-generated costates is shown in Fig.~\ref{fig:costates-ML1}. It apparent from this plot that the costates satisfy \eqref{eq:ML-V&V1}.
%
%

\subsection{A Problem With Discontinuous and Dirac-Delta Functions for Covector Trajectories}
Consider the following deceptively simple optimal control problem\cite{ross-book,brysonHo}:
\begin{eqnarray}
&\bx := (x, v) \in \real{2} \quad \bu := u \in \Real & \nonumber \\[1em] 
(B) &\left\{
\begin{array}{lrl}
\textsf{Minimize } & J[\bx(\cdot), \bu(\cdot)] =&  \displaystyle\frac{1}{2}\int_{t_0}^{t_f} u^2(t)\, dt \\
\textsf{Subject to}& \dot x =& v \\
& \dot v =& u \\
& (x_0, v_0, t_0)  = & (0, 1, 0) \qquad \\
& (x_f, v_f, t_f) =& (0, -1, 1)\\
& x(t) \le & \ell = 0.1
\end{array} \right.&\label{eq:Prob-Breakwell}
\end{eqnarray}
According to \cite{brysonHo}, this problem (with a generic $\ell$) was originally designed by Breakwell to illustrate the various nuances of state variable inequality constraints.  The solution provided in \cite{brysonHo} is based on the original Gamkrelidze-version of Pontryagin's Principle.  The solution provided in \cite{ross-book} is based on the Dubovitskii-Milyutin generalization\cite{dmitruk-JOTA-2017} of Pontryagin's Principle. The basic thrust of this generalization is given by \eqref{eq:HAMVET}. Using \eqref{eq:HAMVET} and \eqref{eq:complementarity}, it can be shown (see Problem~4.4.3 in \cite{ross-book}) that the co-position trajectory, $t \mapsto \lambda_x(t)$ is given by,
\begin{equation}\label{eq:ode=Dirac-delta}
\dot\lambda_x(t) = - \mu(t) = - 22.22\,\delta_{Dirac}(t-t_i), i = 1, 2, \quad t_1 = 0.3, \ t_2 =0.7
\end{equation}
where, $\mu(t)$ is the path covector associated with the state constraint, $x(t) \le 0.1$. As indicated in \eqref{eq:ode=Dirac-delta}, $\mu(t)$ is given explicitly in terms of the shifted Dirac-delta functions, $\delta_{Dirac}(t-t_i)$.  Furthermore, it can be shown\cite{ross-book} that $t \mapsto \lambda_x(t)$ is piecewise constant with downward jumps at $t_1$ and $t_2$.

\begin{figure}[h!]
      \centering
\begin{subfigmatrix}{2} 
\subfigure[Discontinuous staircase costate trajectory $t \mapsto \lambda_x(t)$ generated by a Birkhoff-theoretic computational method]{\includegraphics{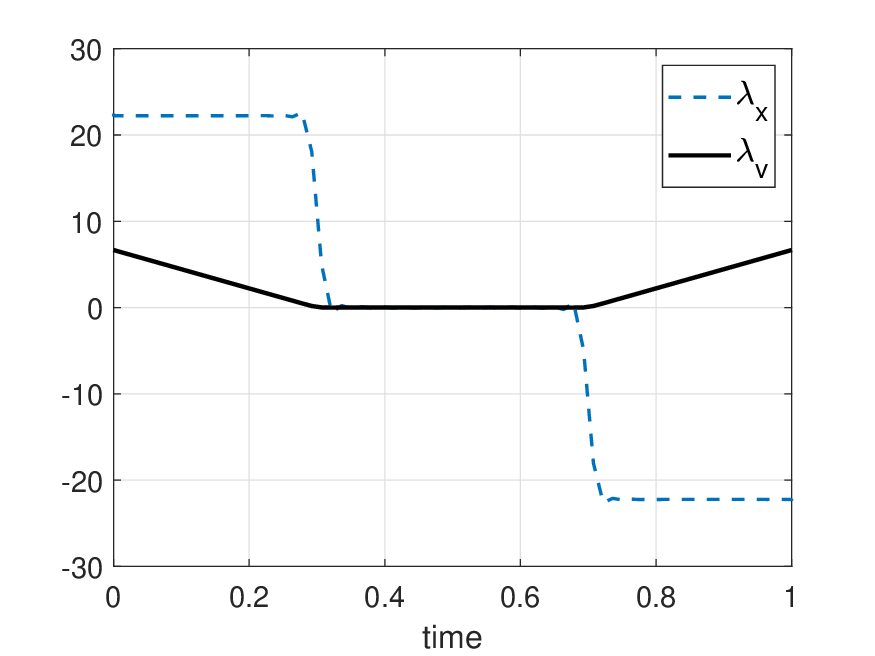}}
\subfigure[A Birkhoff solution to an atomic-measurable covector trajectory, $t \mapsto \mu(t)$, containing two Dirac-delta functions.]{\includegraphics{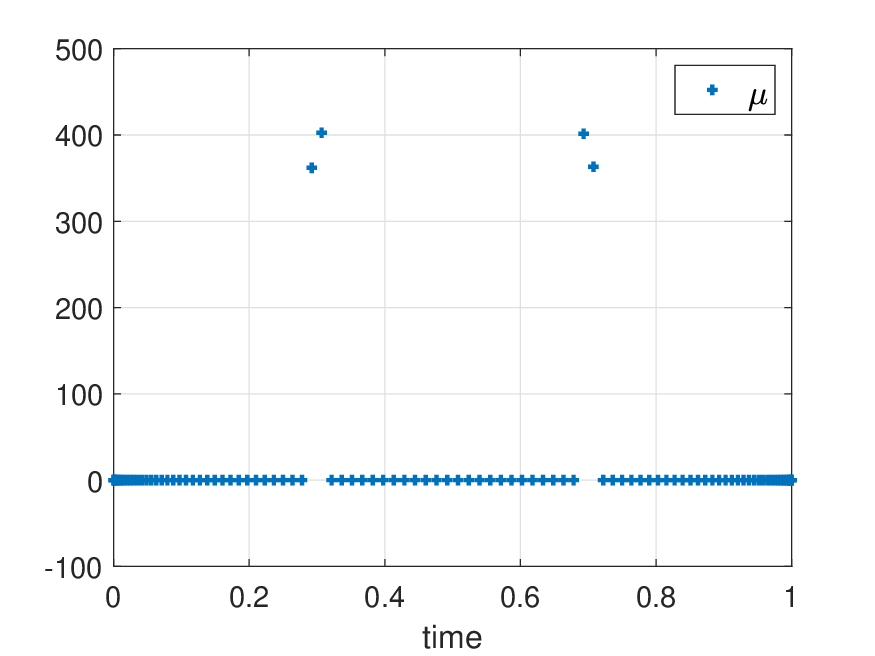}}
\end{subfigmatrix}
\caption{\textsf{Birkhoff-produced discontinuous and and Dirac-delta-type covector trajectories for Problem~$(B)$. }}\label{fig:Breakwell=steps+deltas}
\end{figure}

The Birkhoff-generated covector trajectories to this problem is shown in Fig.~\ref{fig:Breakwell=steps+deltas} for $N =100$.  It is apparent from Fig.~\ref{fig:Breakwell=steps+deltas} that the Birkhoff-computed solution is completely consistent with the theoretical conditions given by \eqref{eq:ode=Dirac-delta}. Note, in particular, the high-accuracy of the atomic measure captured by the Birkhoff solution, $t_k \mapsto \mu(t_k), \ k = 0, \ldots, 100$: just two points around $t_1 = 0.3$ and $t_2 = 0.7$ are off of the zero line in Fig.~\ref{fig:Breakwell=steps+deltas}.b).  These plots demonstrates that a Birkhoff-theoretic method is able to well-approximate Dirac-delta functions and produce jump discontinuities in the costates.
%
%

%


\subsection{A Highly Oscillatory Trajectory Optimization Problem Requiring a Dense Grid}
\begin{figure}[h!]
      \centering
      {\parbox{0.9\textwidth}{
      \centering
      {\includegraphics[width = 0.6\textwidth]{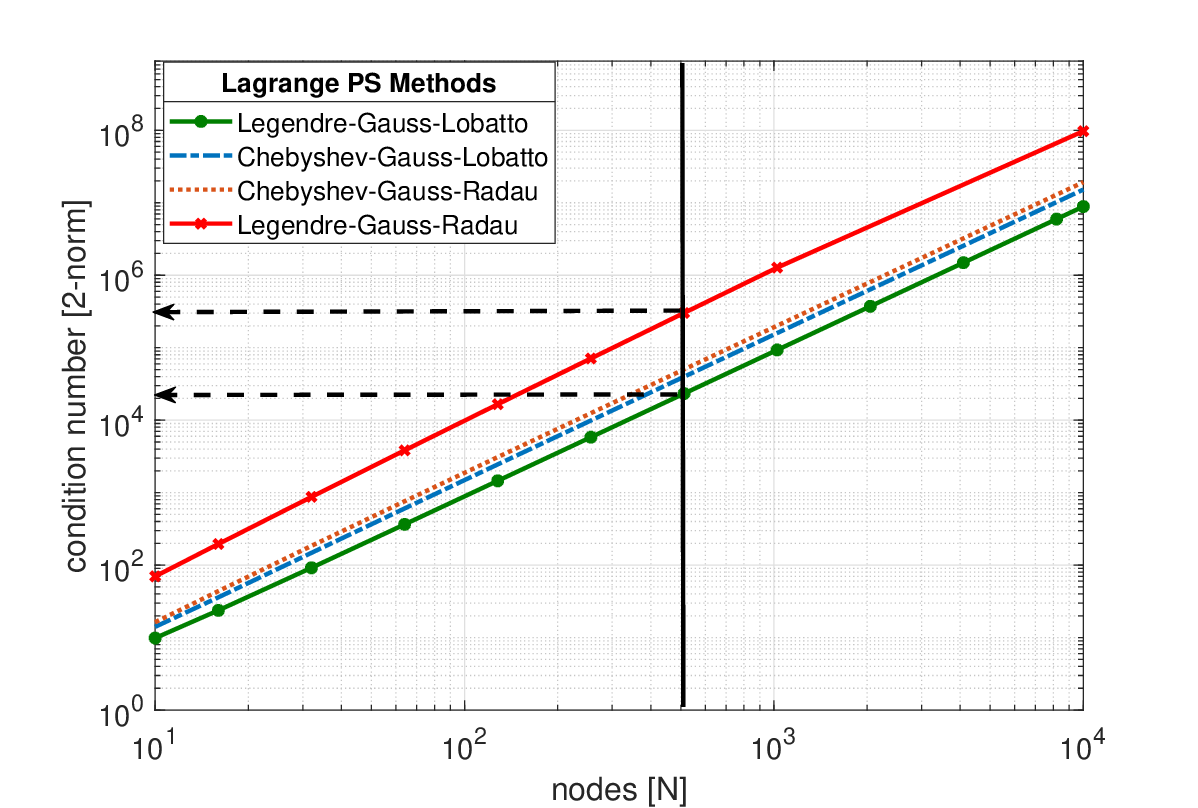}}
      \caption{\textsf{Condition numbers for Lagrange PS methods (cf.~Fig.~1 in ~\cite{newBirk-part-I})}}\label{fig:condNumsLags}
      }
      }
\end{figure}

\begin{figure}[h!]
      \centering
\begin{subfigmatrix}{2} 
\subfigure[Condition numbers for the primal linear system associated with Birkhoff PS methods]{\includegraphics{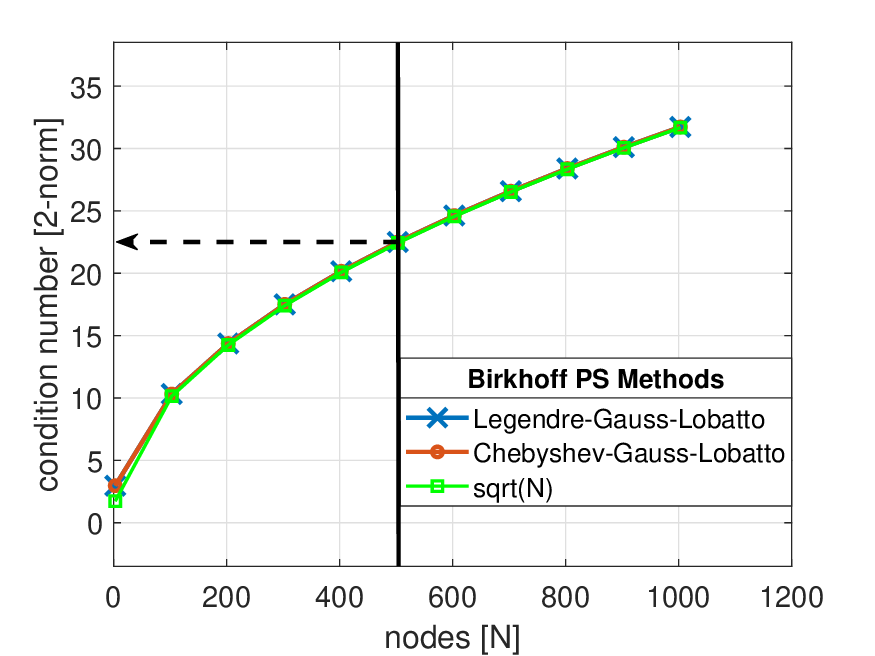}}
\subfigure[Condition numbers for the sub-blocks given by $(I^N, -\bB^a)$]{\includegraphics{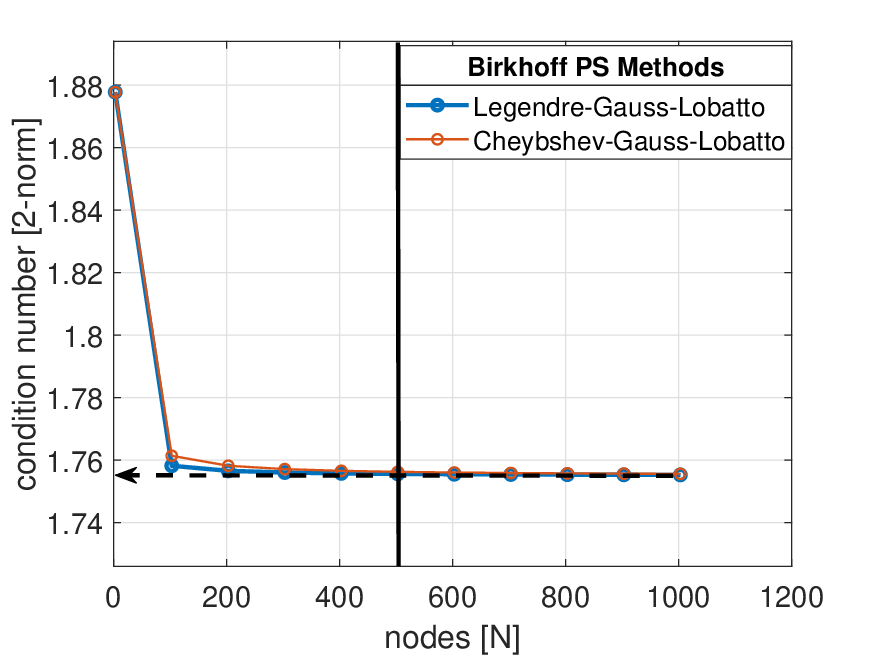}}
\end{subfigmatrix}
\caption{\textsf{Condition numbers of Birkhoff PS methods.}}\label{fig:condNumsBirks}
\end{figure}

A challenge low-thrust circle-to-circle orbit transfer problem is given by\cite{fastmesh},
%
\begin{eqnarray}
& \bx = (r, \theta, v_r, v_t) \in \real 4, \quad \bu = \alpha \in \Real & \nonumber \\ [1em]
(O_\text{Xfer}) &\left\{
\begin{array}{lrl}
\textsf{Minimize } & J[\bx(\cdot), \bu(\cdot), t_f] =& t_f
\\
\textsf{Subject to}
& \dot r =& v_r \\
& \dot\theta =&   v_t/r\\
&\dot v_r =& v_t^2/r - 1/r^2 + 5 \times 10^{-4}\sin\alpha  \\
&\dot v_t =& -v_r v_t/r + 5 \times 10^{-4}\cos\alpha \\
&t_0 = & 0 \\
&\left(r_0, \theta_0, v_{r_0}, v_{t_0}\right) =& \left(1, 0, 0, 1\right)\\
&\left(r_f, v_{r_f}, v_{t_f}\right) =& \left(6, 0, \sqrt{1/6}\right)
\end{array} \right. & \label{eq:prob-OXfer}
\end{eqnarray}
%
where, the numerical values of the various quantities are representative of a LEO-to-GEO maneuver.
What makes this particular problem challenging is that the optimal steering $t \mapsto \alpha(t)$ is highly oscillatory\cite{fastmesh}; hence, the grid must be sufficiently dense to avoid aliasing\cite{bellman-low-t,trefethen-2000,fornberg}.  If one were to use a Lagrange PS method with, say, $N \sim 500$ nodes, it is evident from Fig.~\ref{fig:condNumsLags} (see also Fig.~1 in \cite{newBirk-part-I}) that the ``best'' condition number (LGL grid) would be about $ 10^4$, while the worst (LGR grid) would be about $ 10^5$.
%
%
In sharp contrast, as shown in Fig.~\ref{fig:condNumsBirks}, the condition number of the linear system associated with a Birkhoff PS method (i.e., $[\bA_a(\pi^N), -\bC^N_a]$; see \eqref{eq:ProbbPNa}) is only about 22 for 500 nodes (varies almost exactly as $\sqrt{N}$ as shown in Fig.~\ref{fig:condNumsBirks}.a)) while that of the sub-block given by $[I^N, -\bB^a]$  flatlines to approximately $1.76$.
%
%
As a matter of further contrast, if one were to limit the condition number to approximately $22$, then the maximum number of allowable nodes for a Lagrange PS method would be about $16$ for an LGL grid and less than $10$ for an LGR grid (see Fig.~\ref{fig:condNumsLags}).  Because the CGL/CGR grids have condition numbers that are closer to an LGL grid, the corresponding maximum number of node points for a Chebyshev grid would roughly be the same as that of an LGL grid.  At such low values of $N$, the differences between a high-order Runge-Kutta-type collocation\cite{conway-GL} and Lagrange PS methods become vanishingly small with the latter providing no significant advantages over the former\cite{conway:survey,PSReview-ARC-2012,fornberg,hnw-ode}

\begin{figure}[h!]
      \centering
\begin{subfigmatrix}{2} 
\subfigure[Discrete Birkhoff and interpolated candidate optimal controls (cf.~\eqref{eq:u=grid-soln} and \eqref{eq:uSharp=cont}) to Problem~$(O_{Xfer})$
]{\includegraphics{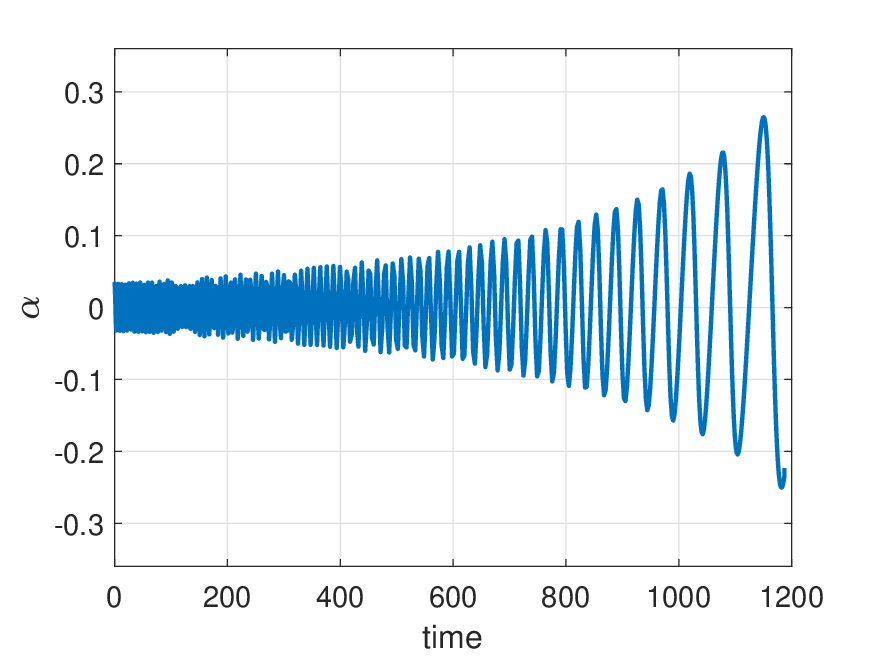}}
\subfigure[Discrete Birkhoff and propagated solutions (cf.~\eqref{eq:x=grid-soln} and \eqref{eq:IVP-V&V}) to Problem~$(O_{Xfer})$ given by \eqref{eq:prob-OXfer}]{\includegraphics{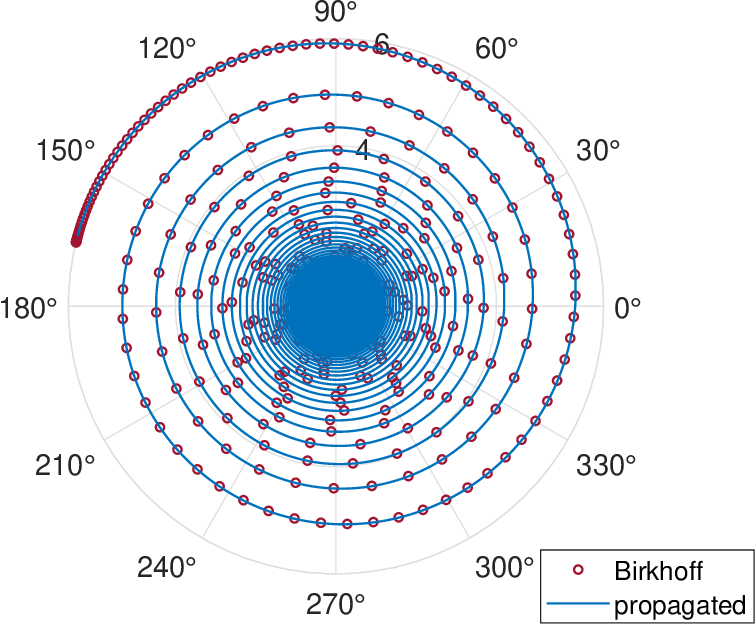}}
\end{subfigmatrix}
\caption{\textsf{A Birkhoff-theoretic primal-feasible solution to Problem~$(O_{Xfer})$ }}\label{fig:spiral+alpha}
\end{figure}

On the basis of the preceding analysis and the data in Fig.~\ref{fig:condNumsBirks}, it is apparent that producing a solution with 500 or more grid points using a Birkhoff method should not wreak havoc in the linear algebra that underpins the trajectory optimization algorithm. Indeed, a candiate primal solution to Problem~$(O_{Xfer})$ for $N=600$ is shown in Fig~\ref{fig:spiral+alpha}.
%
%
It is apparent that the candidate optimal control is highly oscillatory which implies (by the sampling theorem\cite{trefethen-2000,boyd}) that a sufficiently large number of grid points are necessary to adequately represent the solution.  A simple method to test the sufficiency condition is to simply solve the IVP (see \eqref{eq:IVP-V&V}) generated by the candidate optimal control.  A solution to this IVP is the propagated solution indicated in Fig.~\ref{fig:spiral+alpha}.  By a cursory inspection of Fig~\ref{fig:spiral+alpha}, it is clear that this solution is indeed primal feasible (see Definition~\ref{def:primal-feasibility} and Remark~\ref{rem:V&V-primal}).

\begin{figure}[h!]
      \centering
      {\parbox{0.9\textwidth}{
      \centering
      {\includegraphics[width = 0.44\textwidth]{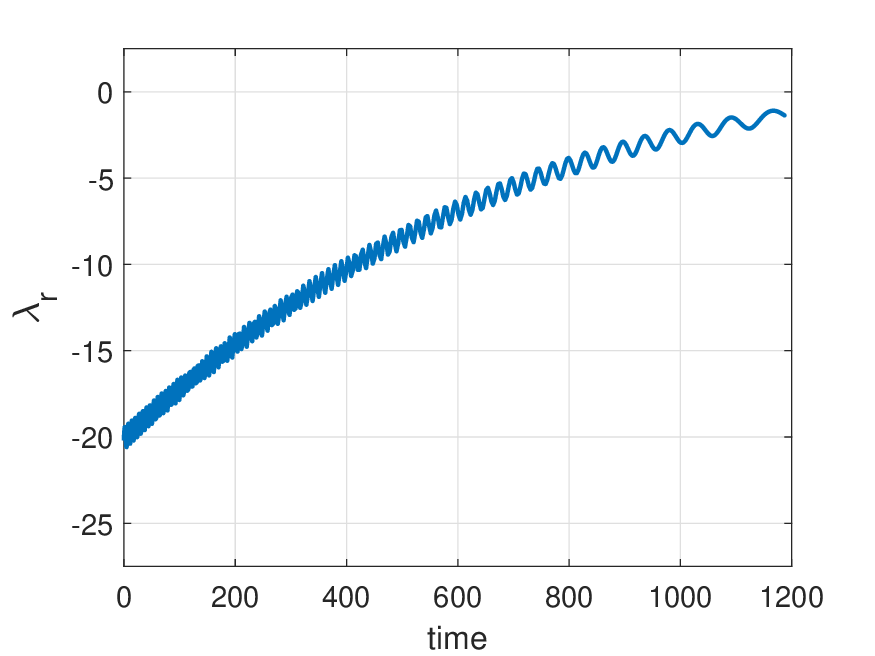}}
      {\includegraphics[width = 0.44\textwidth]{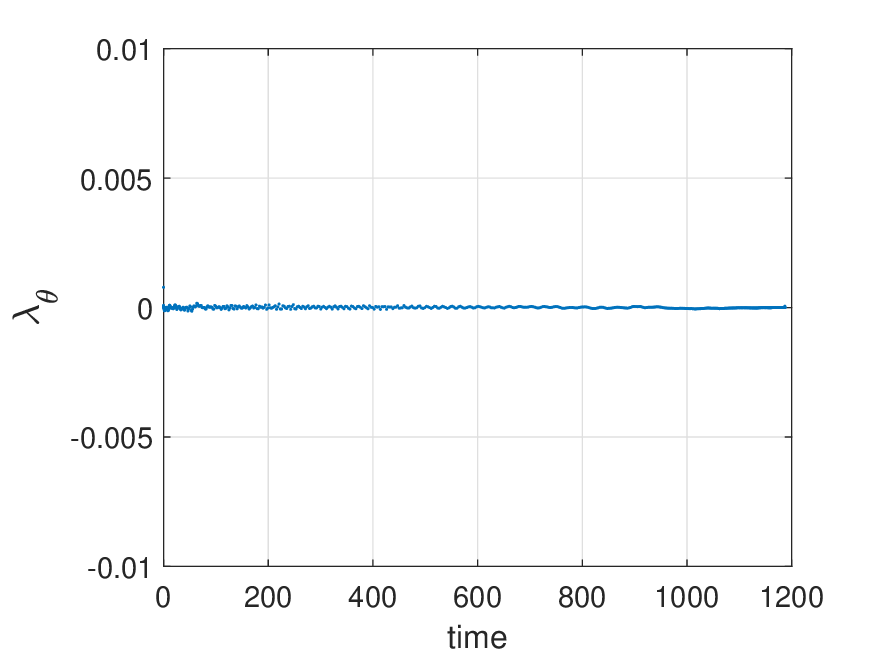}}
      {\includegraphics[width = 0.44\textwidth]{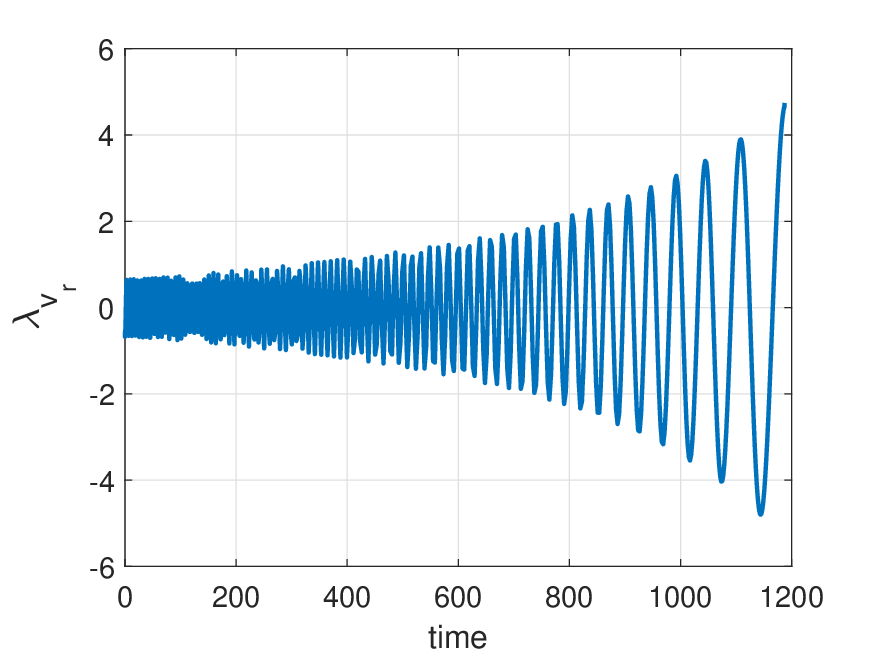}}
      {\includegraphics[width = 0.44\textwidth]{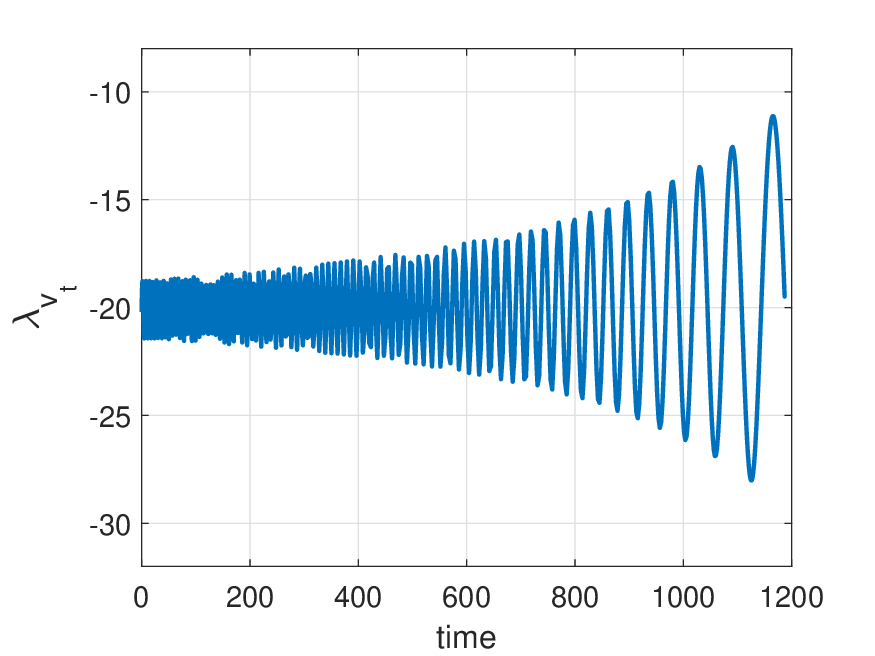}}
      \caption{\textsf{Highly oscillatory costate trajectories for Problem~$(O_{Xfer})$ balanced by a scale factor of $100$ in accordance with the procedure of \cite{scaling}}.}\label{fig:OXfer-duals}
      }
      }
\end{figure}

Beyond the issues pertaining to dense grids and condition numbers, it turns out that part of the challenges in generating an optimal low-thrust solution is that the problem is not well balanced\cite{scaling}.  That is, even under canonical scaling (as already performed in \eqref{eq:prob-OXfer}) it turns out that the costates have very large values that range from zero to $2000$.  In principle, it is possible to generate a well-balanced computational problem by choosing designer units\cite{ross-book} as detailed in \cite{scaling}; however, as shown in Fig.~\ref{fig:OXfer-duals}, a reasonably well-balanced Problem~$(O_{Xfer})$ was feasible by simply scaling down the cost function by a factor of $100$.
%
%
Thus, for example, the values of $\lambda_r$ in canonical-consistent units (see \cite{ross-book,scaling}) has a range of about $2,000$ time units per distance units. By scaling down the cost function (time units) by a factor of $100$, the range of values of $\lambda_r$ was brought down roughly consistent with that of $r$ (of $6$ units, per \eqref{eq:prob-OXfer}).  As apparent from Fig.~\ref{fig:OXfer-duals}, the same scale factor generated a roughly balanced problem with respect to the other primal-dual pairs, namely, $(\theta, \lambda_\theta)$, $(v_r, \lambda_{v_r})$ and $(v_t, \lambda_{v_t})$.  Note also that, except for $\lambda_\theta$, the costates are highly oscillatory.  Taken all together, i.e., grid density, oscillatory functions, condition numbers, scaling and balancing, it is clear why a low-thrust trajectory optimization problem is challenging.

\begin{figure}[h!]
      \centering
      {\parbox{0.9\textwidth}{
      \centering
      {\includegraphics[width = 0.44\textwidth]{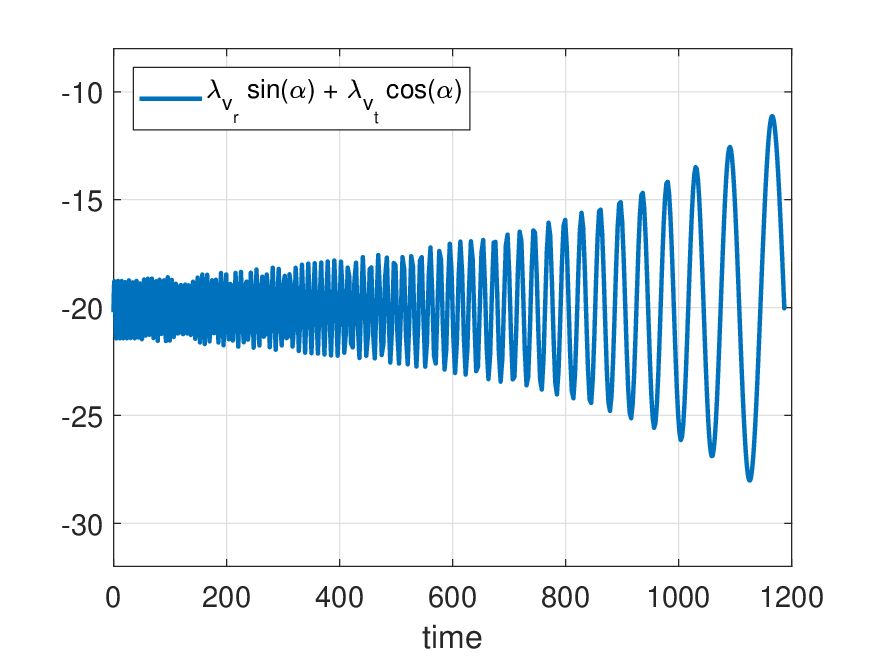}}
      {\includegraphics[width = 0.44\textwidth]{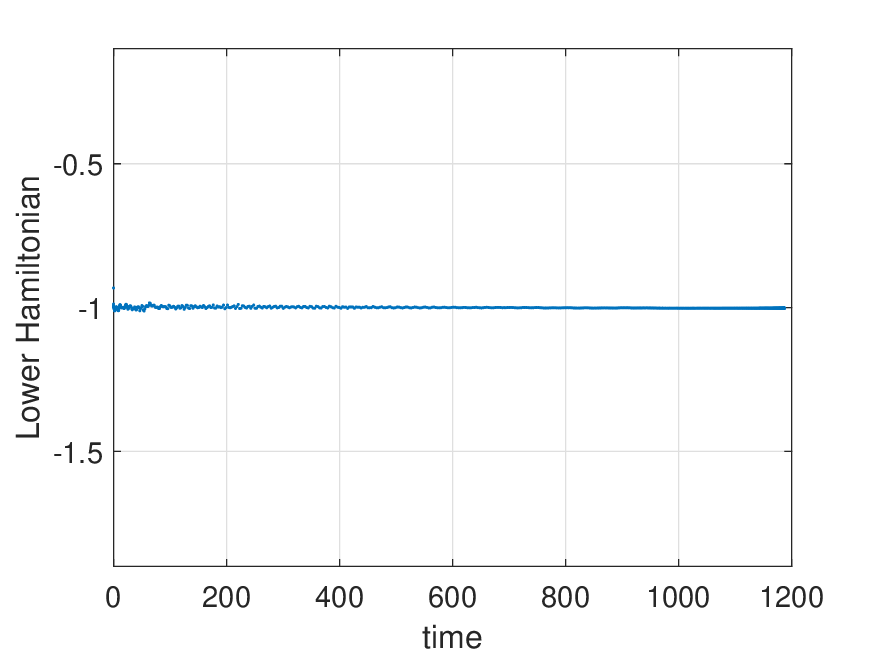}}
      \caption{\textsf{Plots of the left-hand-side of \eqref{eq:Hconvexity} and \eqref{eq:H=-1}}}\label{fig:OXfer-dual-VV}
      }
      }
\end{figure}

Applying \eqref{eq:HAMVET} to Problem~$(O_{Xfer})$, it is straightforward to show that an optimal trajectory must satisfy the following necessary conditions:
\begin{subequations}\label{eq:OXfer-dual-VV}
\begin{align}
\lambda_\theta(t) &= 0 \label{eq:lam_theta=0}\\
\lambda_{v_r}(t)\sin\alpha(t) +\lambda_{v_t}(t)\cos\alpha(t) &\le 0 \label{eq:Hconvexity}\\
\mathcal{H}[@ t] &= - 1 \label{eq:H=-1}
\end{align}
\end{subequations}
That \eqref{eq:lam_theta=0} is satisfied is apparent from Fig.~\ref{fig:OXfer-duals}.  The satisfaction of \eqref{eq:Hconvexity} and \eqref{eq:H=-1} over the primal-feasible solution of Fig.~\ref{fig:spiral+alpha} is shown in Fig.~\ref{fig:OXfer-dual-VV}.
%
%
These dual optimality tests indicate that the primal feasible solution satisfies all of the verifiable optimality conditions resulting from an application of Pontryagin's principle.  As a matter of completion, we briefly note that these verification tests were made possible by the results of Part I of this paper\cite{newBirk-part-I}.

\subsection{A Brief Discussion on the Totality of Results}
At this juncture, it is instructive to provide a broader perspective on all of the numerical results presented. Beyond the specifics illustrated in each of the preceding example problems and solutions, the constant theme in all cases are independent demonstrations of feasibility and the satisfaction of the necessary conditions of optimality via the methods described in Section~\ref{sec:IV&V}.  It is apparent by the presentation of the results in all these cases that the CMP (derived in \cite{newBirk-part-I}) plays a critically important but well-hidden role.  In fact, this is precisely the point of the CMP: to be not distractive in analyzing a particular problem as indicated in the left ``analysis box'' of Fig.~3 in \cite{newBirk-part-I}.

Consider, for the purposes of argument, a direct method that does not satisfy the CMP.  The most that can be said of a solution resulting from such a method is primal feasibility (provided it is independently validated via the methods described in Section~\ref{sec:IV&V}).  It is fairly straightforward to show, via simple examples and flight results, that the gap between feasibility and optimality can easily be 50\% in propellant\cite{ross-howto-2004} (100\% in the case of the Zero-Propellant Maneuver\cite{zpm:NASA-report}) or 70\% in maneuver time\cite{fleming-jgcd-rigid-body} (100\% in the case of ``FastMan,'' a time-optimal maneuver\cite{fastman-2019} implemented onboard the Lunar Reconnaissance Orbiter\cite{LRO:UnivToday2021}).  In other words, the notion that ``feasibility is sufficient for most engineering problems'' is highly misguided.  We also emphasize that finding optimal solutions has intrinsic value even if the optimal solution is never implemented in flight: it informs engineers and operators the price of their suboptimal/nonoptimal decisions. We finally note that a direct method (that does not satisfy the CMP) can also potentially generate false infeasibilities that might produce unsafe or dangerous human flight operations\cite{bollino:PhD-2006}.

\section{Conclusions}
In this part II of a two-part paper, we showed that the universal Birkhoff theory advanced in part I holds for a family of Legendre and Chebyshev grid points, be bit Lobatto, Radau or Gauss. Recent advancements in spectral methods prove that all of these grid points can be generated at an $\mathcal{O}(1)$ computational speed. Nonetheless, the collection of Chebyshev grid points are the easiest to implement because it requires just about a single line of code based on the cosine function.  Among the set of Birkhoff matrices generated for Lobatto, Radau and Gauss grid points, the Lobatto set has the particular property that the last row of the $a$-form of the Birkhoff matrix is identically equal to the Birkhoff quadrature weights.  This convenience implies that the Birkhoff method over a Chebyshev-Gauss-Lobatto grid
is the easiest to implement.  The other five grid selections appear to be equally effective but require additional effort.

The computation of all Birkhoff matrices can be performed rapidly, efficiently and stably over thousands of grid points without suffering round-off errors that plague Lagrange pseudospectral methods; i.e., methods based on differentiation matrices.  All Birkhoff matrices over the Legendre and Chebyshev grids flatten the rapid growth in the condition numbers from $\mathcal{O}(N^2)$ to $\mathcal{O}(1)$.  As proved in this paper, the Birkhoff methods generated over six well-known grid points satisfy the covector mapping principle; i.e., a ``direct'' Birkhoff method is equivalent to an ``indirect'' method.  This equivalence is exploited to reuse previously-developed spectral algorithms that rely on covector mapping theorems. Furthermore, as proved in part I of this paper, the Birkhoff-specific primal-dual system can be isolated to a linear system even for nonlinear problems.  Because this property does not hold for Lagrange pseudospectral methods, the resulting spectral algorithm can be further customized to take advantage of this linear-nonlinear split.

It is clear from the details provided in this paper that a fast and accurate implementation of a Birkhoff method requires an integration of several disparate concepts in mathematical analysis and computational theory.  This daunting task can be alleviated by selecting a Chebyshev-Gauss-Lobatto grid.  This is because a Birkhoff method over this selection is the easiest to implement through the use of readily available off-the-shelf software. As illustrated in Section~\ref{sec:examples} of this paper, a Birkhoff method over a Chebyshev-Gauss-Lobatto grid is fully capable of solving practical problems of interest in aerospace and mechanical engineering even when the primal-dual solution involves discontinuities and Dirac-delta functions.

\section{Acknowledgment}

Partial funding for this research, provided by the Defense Advanced Research Project Agency, is gratefully acknowledged.




\end{document}